\documentclass{article}

\usepackage{graphics,amsmath,amsfonts,amssymb,graphicx,url,color,amsthm,cite}
\usepackage[font=small]{caption}
\usepackage[colorlinks=false]{hyperref}

\usepackage{calc}
\newlength{\notewidth}
\newtheorem{lemma}{Lemma}[section]

\newtheorem{theorem}[lemma]{Theorem}

\theoremstyle{definition}

\newtheorem{remark}[lemma]{Remark}

\newcommand{\eps}{{\varepsilon}}

\newcommand{\proofend}{\hfill$\Box$\bigskip}

\newcommand{\C}{{\mathbb C}}

\newcommand{\R}{{\mathbb R}}

\renewcommand{\v}{\mathbf{v}}

\newcommand{\Rot}{\R^{1,2}}
\newcommand{\st}{\, | \,}
\newcommand{\be}{\begin{equation}}
\newcommand{\ee}{\end{equation}}
\newcommand{\F}{\mathcal{F}}

\begin{document}

\title{Variations on the Tait-Kneser theorem}

\author{Gil Bor\footnote{
CIMAT, A.P. 402, Guanajuato, Gto. 36000, Mexico; 
gil@cimat.mx
}
\and
Connor Jackman\footnote{
CIMAT, A.P. 402, Guanajuato, Gto. 36000, Mexico; 
connor.jackman@cimat.mx
}
\and
Serge Tabachnikov\footnote{
Department of Mathematics,
Penn State University, 
University Park, PA 16802;
tabachni@math.psu.edu}
}

\date{}
\maketitle

%\note{added Connor to authors}

\begin{abstract} The Tait-Kneser theorem, first demonstrated by Peter G. Tait in 1896, states that the osculating circles along a plane curve with monotone non-vanishing curvature are pairwise disjoint and nested. This note contains a proof of this theorem using the Lorentzian geometry of the space of circles. We show how a similar proof applies to two variations on the theorem, concerning the osculating Hooke and Kepler conics along a plane curve.  We also prove a version of the 4-vertex theorem for Kepler conics.
\end{abstract}

 \section{ The Euclidean case revisited} \label{sect:circles}

A smooth plane curve with non-vanishing curvature has at every point an {\em osculating circle} which is tangent to the curve at this point and shares the curvature with it. That is, the osculating circles are 2nd order tangent to the curve at every point. At some points the osculating circle may be  tangent to higher order. Such points are called {\em vertices},  and these are the critical points of the curvature.

\begin{figure}[ht]
\centering
\includegraphics[width=.9\textwidth]{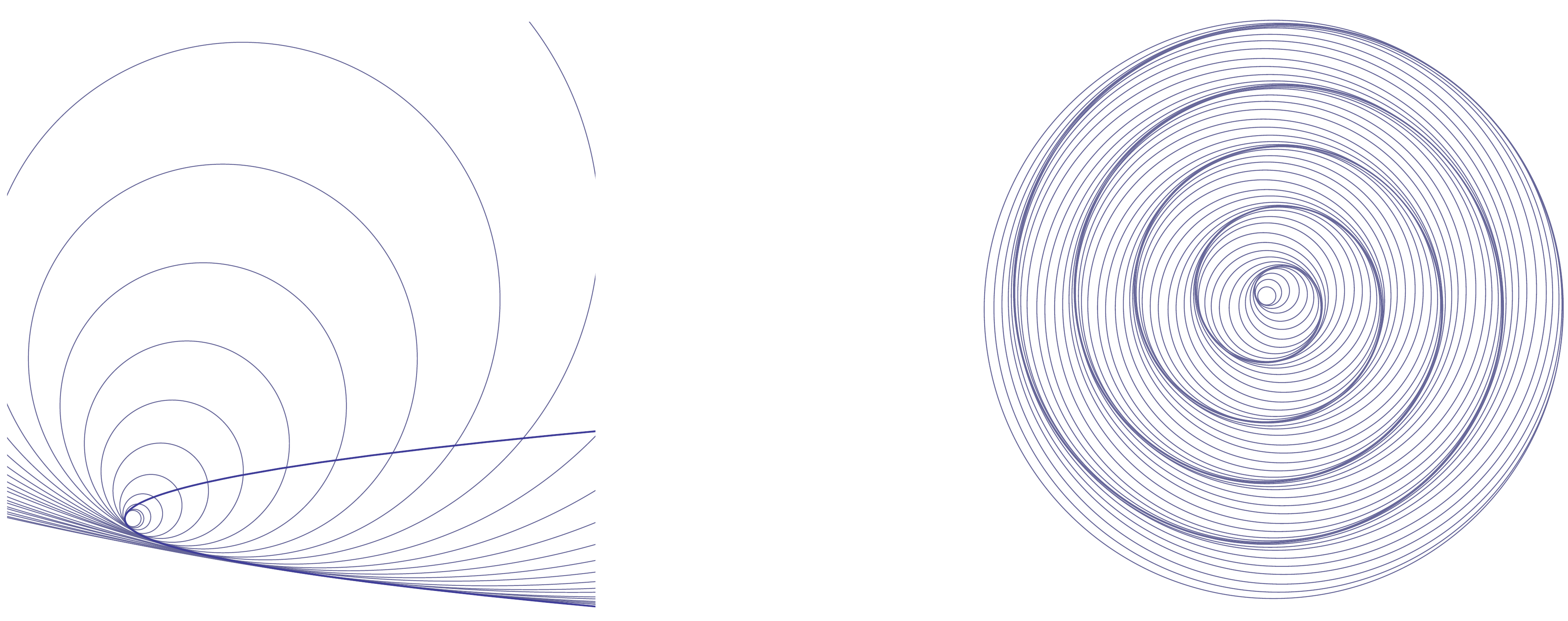}
\caption{Osculating nested circles along a curve with monotone non-vanishing curvature: the lower part of a parabola (left) and an Archimedean spiral (right).}
\label{foliation}
\end{figure}

\begin{theorem}[Tait-Kneser]\label{thm:TK}
The osculating circles of a vertex-free plane curve with non-vanishing curvature are disjoint and pairwise nested, see Figure \ref{foliation}.
\end{theorem}

%{\color{red} Need to redo this picture that was ``borrowed" from another article.}

This theorem is more than a century old \cite{T}; it has numerous variations and ramifications, see the survey \cite{GTT}. 

The osculating circles of a curve with monotone curvature form a foliation of the annulus bounded by the osculating circles at the end points of the the curve. 
%
%\note{added "of the annulus  etc", after the refree asked "a foliation of what?"}
%
We leave it to the reader to mull over the seemingly paradoxical property of this foliation: the curve is tangent to the leaves at every point, but it is not contained in a single leaf (isn't it similar to a non-constant function having everywhere vanishing derivative?!)

The original Tait's proof was very short, just two paragraphs long; it made use of the notions of  evolute and involute of a curve. We present a different proof of Theorem \ref{thm:TK} that will also work in the variations  to follow. 

A circle in  $\R^2$   is given by an equation of the form $(x-a)^2+(y-b)^2=r^2.$ Denote  by $\Rot$  the three dimensional  pseudo-Euclidean space  with coordinates $a,b,r$ equipped with the indefinite quadratic form  $|(a,b,r)|^2:=-a^2-b^2+r^2$ (we use this notation even though the latter may happen to be  a negative number!). The space of  circles  in $\R^2$ is  parametrized by the upper half space $\Rot_+:=\{(a,b,r)\st r>0\}$.  

The null cone with vertex  $\v_0\in\Rot$
is the set of points $\v\in\Rot$ such that $|\v-\v_0|^2=0$, and its interior consist of points $\v$ with $|\v-\v_0|^2>0$. We next describe when two circles are {\em nested}, that is,  the interior of one of them is contained in that of the other.

\begin{lemma} \label{lm:nest}
Given two  circles and the corresponding points $\v_1, \v_2\in\Rot_+$,  the circles are nested  if and only if $|\v_1-\v_2|^2\geq0$, with equality when the circles are nested and tangent.   
\end{lemma}

\begin{figure}[ht]
\centering
\includegraphics[width=\textwidth]{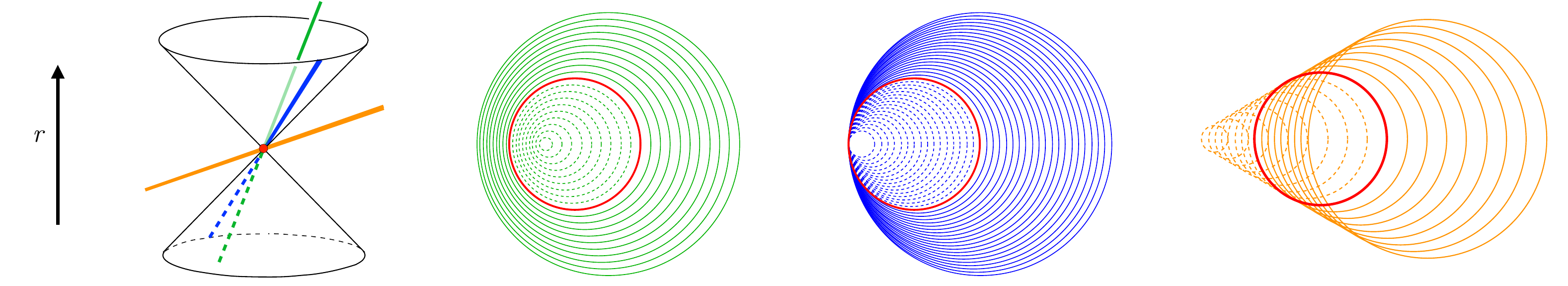}
\caption{`Lines'  of circles. A timelike line (green): nested disjoint circles. A null line (blue):  nested circles tangent at a point. A spacelike  line (orange): intersecting circles, tangent to a pair of lines. }
\label{fig:circles}
\end{figure}

That is, the circles are nested when one of the corresponding points in $\Rot$   lies in  the light  cone   whose vertex is the  other point.  See Figure \ref{fig:circles}. 
 
\proof
Let the radii of the circles be  $R\geq  r$, and the distance between their centers be $d$. 
The nesting condition is $d+r\leq R$,  or  $|\v_1-\v_2|^2=-d^2+ (R-r)^2\geq 0,$ with equality if and only if $ d+r= R$, that is, the circles are tangent. 
\proofend

%{\color{red} Need to redo this picture that was ``borrowed" from another article.}

Let $\gamma$ be a plane curve with non-vanishing curvature and let $\Gamma \subset \Rot$ be the curve of osculating circles of $\gamma$. The vertices of $\gamma$ correspond to singular points of $\Gamma$, and  $\Gamma$ is regular if $\gamma$ is vertex free (see Equation \eqref{eq:reg} in the proof of the next Lemma). 
%
%\note{added the sentence in parenthesis as the referee asked to justify the statement}

A curve in $\Rot$ is {\it null} if its tangent vector is tangent to the null cone at every point. 

\begin{lemma} \label{lm:null}
The curve $\Gamma$ is a null curve.
\end{lemma}

\proof
Let $\gamma(t)=(x(t), y(t))$ be an arc length parameterization of $\gamma$. Then $\kappa=x'y''-y'x''$ is the curvature of $\gamma$, and the  osculating circle  at a point $(x,y)$  of $\gamma$ is given by the  equation  $(X-a)^2+(Y-b)^2=r^2$ with  
$$
(a,b)=(x,y)+r(-y',x'),\ r= {1\over \kappa}.
$$
Since $\Gamma(t)=(a(t), b(t), r(t))$, one has
\begin{align}
\begin{split}\label{eq:reg}
\Gamma'&=-{\kappa'\over\kappa^2}(-y',x',1)+(x',y',0)+{1\over \kappa}(-y'',x'',0) =\\
&= -{\kappa'\over\kappa^2}(-y',x',1),
\end{split}
\end{align}
the last equality due to equation $(x'',y'')=\kappa(-y', x')$. Since $(x')^2+(y')^2=1$, one has $|\Gamma'|^2=0$, as claimed.
\proofend

Here is the intuition behind it with a `hand-waving' proof. An osculating circle $C$ passes through three ``consecutive" points of the curve, say $\gamma(t-\eps), \gamma(t), \gamma(t+\eps)$.  The ``next" osculating circle shares two of these points,
$\gamma(t), \gamma(t+\eps)$, with $C$. In the limit $\eps \to 0$, this implies that the curve $\Gamma$ is tangent to the cone whose vertex is the circle $C$ and consists of the circles tangent to it.    

The last ingredient of the proof of the Tait-Kneser theorem is the next lemma.

\begin{lemma} \label{lm:ineq}
A  regular null curve  $\Gamma:[t_0,t_1]\to \Rot$ satisfies   $|\Gamma(t_1)-\Gamma(t_0)|^2\geq 0$, with equality  if and only if $\Gamma$ is the  null line segment connecting its endpoints.
\end{lemma}

\proof
Let $\Gamma(t)=(a(t), b(t),r(t))$ and $\bar\Gamma(t)=(a(t),b(t))$ be the projection on the horizontal plane. 
The nullity condition on $\Gamma$, $(r')^2=(a')^2+(b')^2,$ implies that the length  of $\bar\Gamma$  is $|r(t_1)-r(t_0)|.$ This length is at least the distance between the endpoints $ |\bar\Gamma(t_1)-\bar\Gamma(t_0)|$,  with equality if and only if $\bar\Gamma$ is the straight line segment. 

It follows that $|\Gamma(t_1)-\Gamma(t_0)|^2=-|\bar\Gamma(t_1)-\bar\Gamma(t_0)|^2+|r(t_1)-r(t_0)|^2\geq 0,$ with equality if and only $\bar\Gamma$ is the line segment. By the nullity condition, this  is equivalent  to $\Gamma$ being a null line segment. 
\proofend

Theorem \ref{thm:TK} follows. Let $\gamma$ be a plane curve with non-vanishing monotone curvature.  Consider two  osculating circles $C_0, C_1$ along $\gamma$ and the regular 
 null curve $\Gamma$   connecting the corresponding points   $\v_0, \v_1\in\Rot$. By the above lemmas,   either $|\v_1-\v_0|^2>0$, in which case $C_0, C_1$ are nested, or $|\v_1-\v_0|^2=0$,  in which case $\Gamma$   is a null  segment connecting $\v_0$ to $\v_1$. The later case corresponds to a family of circles tangent at a point, which cannot occur as the family of osculating circles to a curve in $\R^2$.
 
 \begin{remark}
 {\rm The Lorentzian geometry of the space of circles and the relation between the osculating circles of a curve with null curves is investigated in \cite{NP}.
  }
 \end{remark}

\section{Centroaffine geometry: Hooke orbits} \label{sect:Hooke}

A smooth plane curve $\gamma$ is {\em star-shaped} with respect to the origin if $[\gamma,\gamma']\neq 0$ (the bracket is the determinant made by a pair of vectors). Such a curve can be parameterized so that $[\gamma(t),\gamma'(t)]=1$. Then 
$\gamma''=-p \gamma$, where the function $p(t)$ is  called the {\em centroaffine curvature} of $\gamma$. For example, the origin centered circle of radius $r$  has $p=1/r^4.$ 

A central conic $ax^2+2bxy+cy^2=1$ has  $p=ac-b^2$, the determinant of the coefficient matrix ${a\ b\choose b\ c}$. Central conics are the trajectories of mass points subject to the Hooke force law: the radial force is proportional to the distance to the origin. If the force is attractive the trajectory is a central  ellipse with  $p>0$, and if it is repulsive the trajectory is a hyperbola with $p<0$. Central conics  play the role of circles in centroaffine geometry. 

The {\em osculating central conic} of a star-shaped curve $\gamma$ with non-vanishing centroaffine curvature is the central conic  tangent to $\gamma$ and sharing its centroaffine curvature at the tangency point. It coincides with $\gamma$ to 2nd order at the point of tangency, and the order is higher if 
$p'=0$ at this point. 

Here is a centroaffine version of the Tait-Kneser theorem.

\begin{figure}[ht]
\centering
\includegraphics[width=.9\textwidth]{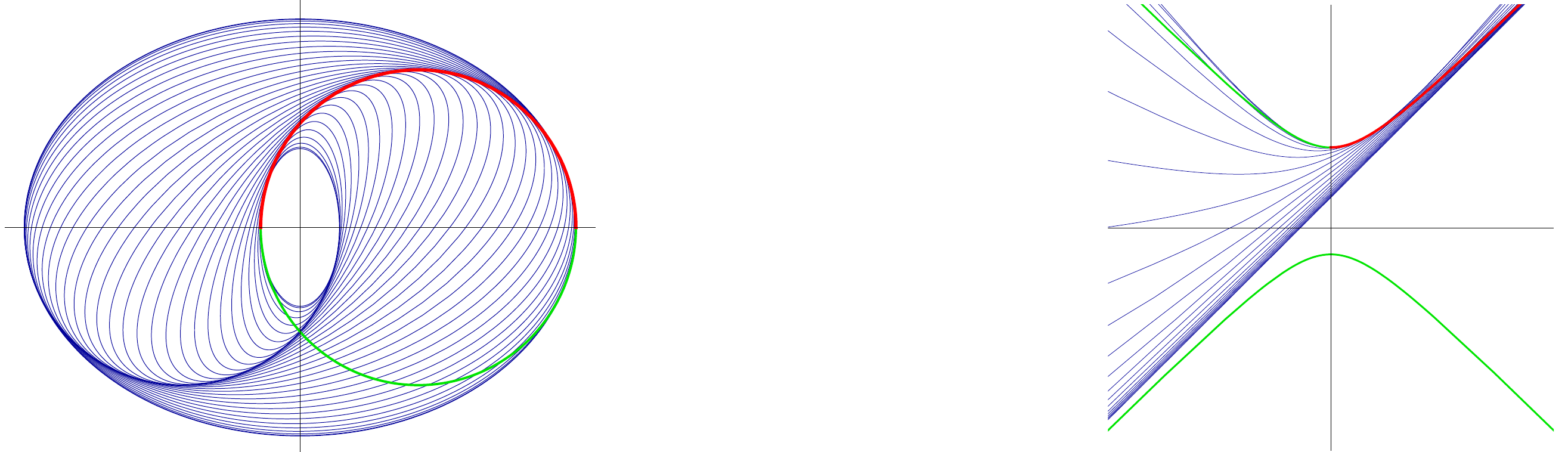}
\caption{Osculating central conics. Left: along the upper part (red) of the circle $(x-.75)^2+y^2=1$ (green). Right: along the upper right part (red) of the hyperbola  $(y-.5)^2-x^2=1$ (green).}
\label{osc_conics}
\end{figure}

 \begin{theorem} \label{thm:TKca}
 The osculating central conics of a star-shaped plane curve with monotone non-vanishing centroaffine curvature are pairwise disjoint and nested, see Figure \ref{osc_conics}.
 \end{theorem}
 
The proof of Theorem \ref{thm:TKca} goes along the same lines as our proof of the Tait-Kneser theorem.

The space of central conics is 3-dimensional with coordinates $(a,b,c)$. This space has a pseudo-Euclidean metric given by the determinant of the quadratic form, $ac-b^2$, and we identify it with $\Rot$.

 An analog of Lemma \ref{lm:nest} holds:

\begin{lemma} \label{lm:nestca}
Let $C_1, C_2$ be two central conics of the same type (ellipses or hyperbolas). If $\det(C_2-C_1)\geq 0$ then they are nested, and  equality implies  they are tangent.   (We denote  a conic  and the quadratic form defining it  by  the same letter.)
\end{lemma}

\proof
We use the well known fact from linear algebra that  two real quadratic forms,  one of which is definite (positive or negative),  can be diagonalized simultaneously. In  dimension 2, a quadratic form  is definite  if and only if its determinant is positive.

Now if  $C_1,C_2$ are ellipses then  the quadratic forms are both positive definite, so, without loss of generality, by the above fact, the ellipses are $x^2+y^2=1$ and  $ax^2+by^2=1$, with $a,b>0.$  The condition $\det(C_2-C_1)\geq 0$ is then $(a-1)(b-1) \geq 0$, which is equivalent to $a,b\geq 1$ or $a,b\leq 1$. In both cases the ellipses are nested, with equality when they are tangent.

\begin{figure}[ht]
\centering
\includegraphics[width=.9\textwidth]{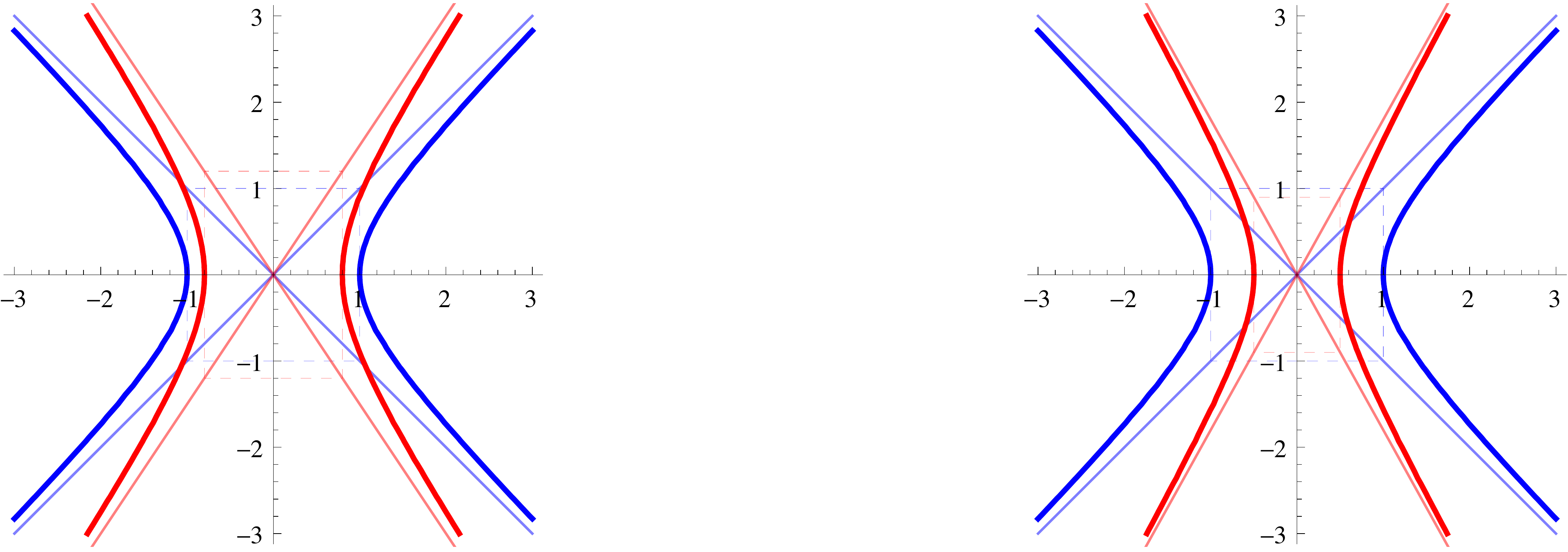}
\caption{The proof of Lemma \ref{lm:nestca}. A pair of nested   hyperbolas is shown in each figure, $x^2-y^2=1$ (blue) and  $(x/a)^2-(y/b)^2=1$ (red).
Left: if the  pair has timelike difference, $\det(C_2-C_1)>0$, then either $0<a<1<b$ or $0<b<1<a$ and the pair is nested. Right: timelike difference is not a necessary nesting condition; if
$a<b<1$  or $1<b<a$ then the difference is spacelike but the hyperbolas are still nested. 
}
\label{fig:hyp}
\end{figure}

If $C_1,C_2$  are hyperbolas, suppose that $\det(C_2-C_1)>0.$  Then $\Delta C:=C_2-C_1$ is a definite quadratic form and, by interchanging the conics if necessary, it is positive. Hence  $\Delta C, C_1$ can be transformed to $ax^2+by^2=1, x^2-y^2=1$ (respectively), with $a,b>0$. Then $C_2$ is $(a+1)x^2-(1-b)y^2=1$, and since it is a hyperbola, we have $0<b<1.$  Renaming  the constants, $C_2$  is $ (x/a)^2-(y/b)^2=1,$  $0<a<1<b.$ It is now easy to see that $C_1$  is nested in $C_2$. See Figure \ref{fig:hyp}.

The case $\det(\Delta C)=0$ is a limit  of  $\det(\Delta C)>0$,  since being ``nested'' is a closed condition. 
\proofend

\begin{remark}Although we do not use it, we note that the converse of Lemma \ref{lm:nestca} holds for ellipses (as follows easily  from the above proof), but not for hyperbolas. For example,  the hyperbolas $x^2-y^2=1$ and  $(x/3)^2-(y/2)^2=1$ are nested, but the respective determinant is negative. We do not dwell on the precise algebraic condition for a pair of quadratic forms to determine a pair of nested hyperbolas. 
\end{remark}

Lemma \ref{lm:null} also has an analog (with the same hand-waving explanation as before).
Again we denote by $\Gamma$ the curve in the space of Hooke conics  that osculate a centroaffine curve $\gamma$.

\begin{lemma} \label{lm:nullca}
The curve $\Gamma$ is a null curve.
\end{lemma}

\proof
Let $\gamma(t)=(x(t), y(t))$ where $xy'-yx'=1$, and let $p(t)$ be the centroaffine curvature. A direct calculation shows that the osculating central conic is given by the equation $aX^2+2bXY+cY^2=1$, with 
$$
a=py^2+(y')^2,\ b=-(p x y+x' y'),\ c=p x^2+(x')^2
$$
(one needs to check that 
$
ax^2+2bxy+cy^2=1,\ (ax+by)x'+(bx+cy)y'=0,\ {\rm and}\ p=ac-b^2).
$ 

Then another calculation shows that 
$$
a'=p'y^2,\ b'=-p'xy,\ c'=p'x^2,
$$
and hence $a'c'-(b')^2=0$, as claimed. 
\proofend

Now Lemma \ref{lm:ineq} applies, both for ellipses and hyperbolas, thereby completing the proof of Theorem \ref{thm:TKca}. 

Let us add that the Tait-Kneser theorem is closely related to another classical result, the 4-vertex theorem that, in its simplest form, states that a plane oval has at least four vertices. As the example in Figure \ref{osc_conics} shows, its analog does not hold in centroaffine geometry: the circle of radius 1 centered at $(0.5,0)$ has only 2 ``vertices", that is, hyper-osculating central conics. 

\begin{remark}
{\rm The Lorentzian geometry of the space of central conics is studied in \cite{Sa}.
}
\end{remark}

\section{Kepler orbits} \label{sect:Kepler}
 
\begin{figure}[ht]
\centering
\includegraphics[width=.9\textwidth]{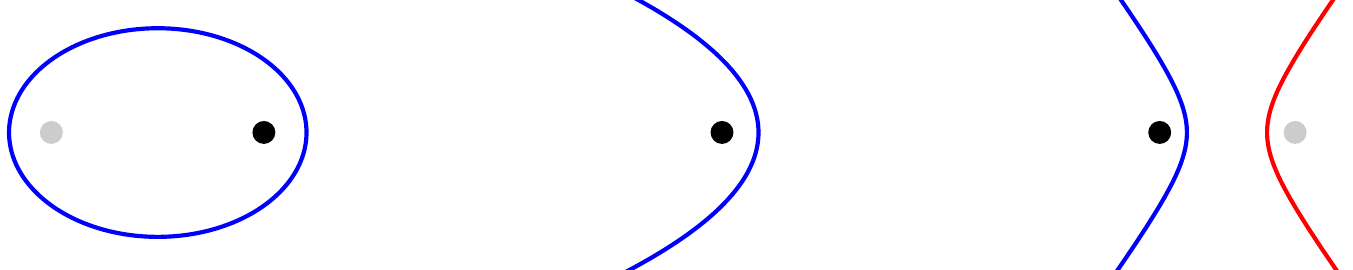}
\caption{Kepler orbits share a focus at the origin (the black dot; the gray dot is the  other focus, not shared). The blue curves (ellipse, parabola, hyperbola) bend around the origin, due to an attractive force. The red curve (hyperbola) bends away from the origin, due to  a  repulsive force.}
\label{fig:kepler_conics}
\end{figure}

A Kepler orbit is a plane conic (ellipse, parabola or hyperbola)  
with a focus at the origin. These are the trajectories of mass points subject to Newton's force law (either attractive or repulsive): the radial force is proportional to the inverse square of the distance to the origin. For attractive force the orbits are ellipses, 
parabolas or hyperbola branches bending around the origin. For repulsive force only hyperbolas appear,  the `other' branches left out by the attractive force hyperbolas. See Figure \ref{fig:kepler_conics}.

%\note{added a reference to Givental's article, after referee request. There should be an older reference, but I dont know any other. }
%
A Kepler conic is the orthogonal projection on the horizontal plane of  the intersection of the cone $x^2+y^2=z^2$ in $\R^3$ with a plane $ax+by+cz=1,\ c> 0$ \cite{Giv}. 
Thus the space of Kepler conics is parametrized by the space  $\Rot_+$ with coordinates $(a,b,c)$ and the quadratic form  $|(a,b,c)|^2=-a^2-b^2+c^2$. The null cone $|(a,b,c)|^2=0$ parametrizes Kepler parabolas, its interior $|(a,b,c)|^2>0$ parametrizes ellipses, and its exterior $|(a,b,c)|^2<0$ parametrizes hyperbolas.

%\note{I think we  need to restrict to curves with positive  centroaffine curvature, since all Kepler orbits have positive curvature. Another option is to throw in also the 2nd branch of hyperbolas (coming from repulsive Coulomb force). But even with this we need to exclude vanishing centroaffine curvature. }
 As before, a smooth star-shaped curve with non-vanishing centroaffine curvature can be 2nd order approximated at every point by its osculating Kepler conic;  this osculating  conic may hyper-osculate, that is, approximate the curve to higher order.

\begin{figure}[ht]
\centering
\includegraphics[width=\textwidth]{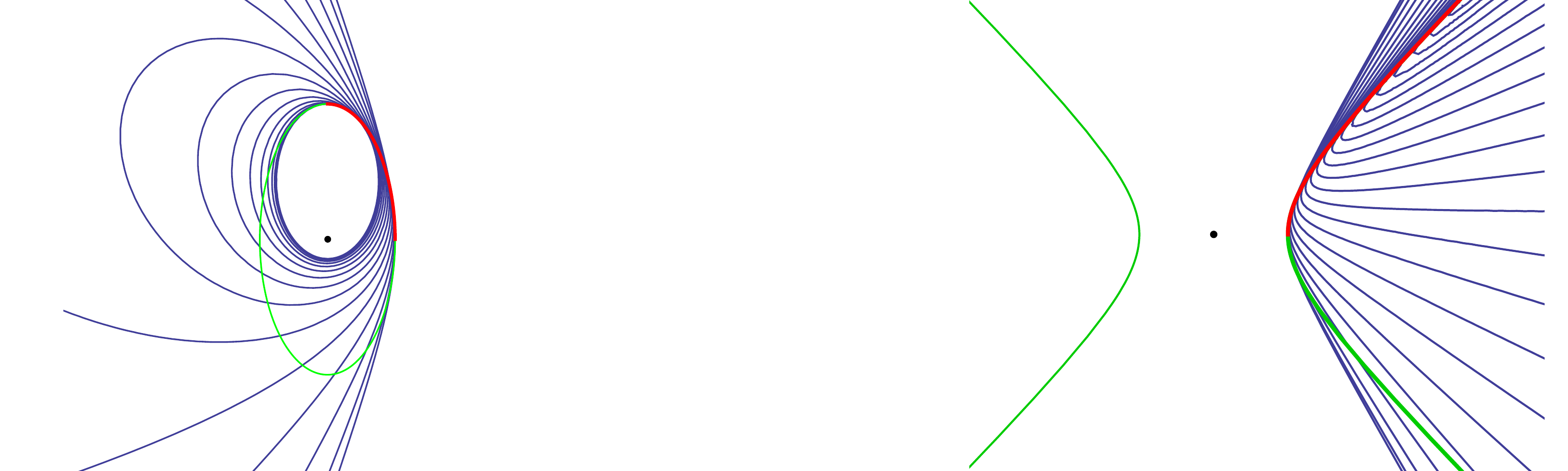}
\caption{Osculating Kepler conics (blue) along a vertex free arc (red) of a Hooke conic (green). Left:  along  a Hooke ellipse. Right: along  a Hooke hyperbola.}
\label{osc_kepl_conics}
\end{figure}

An analog of the Tait-Kneser theorem holds:

%\note{dont we miss a ``non-vanishing curvature" type condition? i.e.,  that the curve can be approximated by osculating Kepler conics. For example, a straight affine line is star-like but shouldnt be covered by this theorem.}
\begin{theorem} \label{thm:TKKe}
Consider a star shaped curve, free from hyper-osculating Kepler conics. Then its osculating Kepler conics are pairwise nested. 
See Figure \ref{osc_kepl_conics}. 
 \end{theorem}
 
 %{\color{red} Would be good to include a figure obtained from Fig. 2 by $z\mapsto z^2$ to illustrate the theorem.}
 
 One can prove Theorem \ref{thm:TKKe} along the same lines as before. We only present an analog of Lemmas \ref{lm:nest} and \ref{lm:nestca}.
 
\begin{lemma} \label{lm:nestKe}
Two Kepler conics, corresponding to points $\v_1,\v_2\in \Rot_+$, are disjoint if and only if $|\v_1-\v_2|^2>0$.
\end{lemma}
  
  (For Kepler conics, being nested and disjoint are equivalent).

\begin{figure}[ht]
\centering
\includegraphics[width=.25\textwidth]{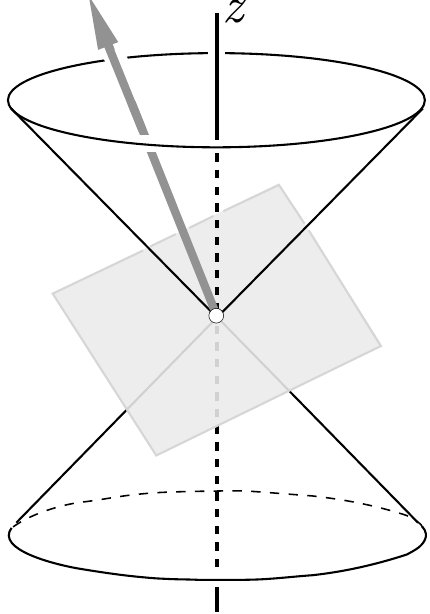}
\caption{Proof of Lemma \ref{lm:nestKe}. The plane $ax+by+cz=0$ intersects the cone $x^2+y^2=z^2$ only at its vertex whenever the normal vector $(a,b,c)$ lies inside the cone.}
\label{cone_plane}
\end{figure}

\proof
The system of equations
$$
a_1x+b_1y+c_1z=1,\ a_2x+b_2y+c_2z=1,\ x^2+y^2=z^2
$$
has no solutions if and only if the system
$$
(a_1-a_2)x+(b_1-b_2)y+(c_1-c_2)z=0,\ x^2+y^2=z^2
$$
has no non-zero solutions,  if and only if the vector $(a_1-a_2, b_1-b_2,c_1-c_2)$ lies in the interior of the cone $x^2+y^2=z^2$,
 as needed. See Figure \ref{cone_plane}. 
\proofend

Unlike Hooke conics, one has a version of the 4-vertex theorem for Kepler conics. Let $\gamma$ be a simple closed star-shaped curve, and let us call a vertex a point at which it is approximated by a Kepler conic to 3rd order (higher than usual). 

\begin{theorem} \label{thm:4v}
The curve $\gamma$ has at least four distinct vertices.
\end{theorem}

\proof
A Kepler conic, in polar coordinates $(\alpha,r)$, is given by the formula
$$
r=\frac{c}{1+e \cos (\alpha+\varphi)},
$$
where $c,e$, and $\varphi$ are constants. Let $\rho=1/r$. Then $\rho(\alpha)$ is a first harmonic, hence $\rho'''+\rho'=0$ for all $\alpha$. This equation characterizes Kepler conics.

The curve $\gamma$ is given by its function $\rho(\alpha)$, and its vertices are precisely the points where $\rho'''+\rho'=0$. It remains to use the fact that such an equation has at least four distinct roots for every periodic function $\rho$ -- a particular case of the Sturm-Hurwitz theorem that asserts that a smooth $2\pi$-periodic function has at least $2n$ roots where $n$ the number of the first harmonic in its Fourier expansion (see, e.g., \cite{OT}, Appendix 1). The case at hand is $n=2$: the function $\rho'''+\rho'$ is free from the constant term and the first harmonics.
\proofend

\section{Odds and ends} \label{sect:odds}

%\note{I split  the 1st sentence  in two and added to the 1st part  references, as asked by the referee.} 
%
Identifying the plane with $\C$, consider the square map $z\mapsto z^2$. This map takes Hooke conics to Kepler conics \cite{Bo, Ar}. It thus provides a direct connection between the results of Sections \ref{sect:Hooke} and \ref{sect:Kepler}. 

A more general statement relates the trajectories of mass points subject to radial forces proportional to powers of distances to the origin; see \cite{Ar} for a modern treatment.

\begin{theorem}[Bohlin-Kasner \cite{Bo,Ka2}] \label{thm:Bo}
Consider two central force laws  in the plane, with the force proportional to $r^a$ and to $r^b$, where $r$ is the distance to the origin. Let $(a+3)(b+3)=4$. Then the map $z\mapsto z^{(a+3)/2}$ takes the trajectories of motion in the first field to those in the second field.
\end{theorem}

The Hooke and Newton attraction laws are  $a=1$ and $b=-2$ (respectively). These cases are distinguished among central force  laws.

\begin{theorem}[Bertrand \cite{Be}] \label{thm:Be}
Assume that all the trajectories of a mass point, subject to a central force that depend on the distance to the origin and whose energy does not exceed a certain limit, are closed. Then the law of attraction is either Hooke's or Newton's. 
\end{theorem}

One notes that the family of circles, of Hooke conics, and of Kepler conics each depends on three parameters -- this is why they can approximate smooth curves to 2nd order.  

More generally, given a field of forces in the plane, the trajectory of a mass point depends on its initial position and velocity, and hence the trajectories form a 3-parameter family. Which 3-parameter families of plane curves are obtained this way? This problem was thoroughly studied by E. Kasner who obtained a complete answer to this question \cite{Ka1,Ka2}.

Other examples of 3-parameter families of curves for which a version of the Tait-Kneser theorem holds are parabolas $y=ax^2+bx+c$, the graphs of quadratic polynomials, and hyperbolas $(x-a)(y-b)=c^2$, the graphs of fractional-linear functions (see \cite{GTT}). Our proof works in these cases as well, with the metrics given by $ds^2=db^2-4(da)(dc)$ and $ds^2= dc^2-(da)(db)$, respectively.

One can describe the graphs of 3-parameter families of functions $y(x)$ by 3rd order differential equations $y'''=F(x,y,y',y'')$. For example, the equation $y'''=0$ describes the vertical parabolas, and the graphs of fractional-linear functions are described by the vanishing of the Schwarzian derivative:
$$
\frac{y'''}{y'}-\frac{3}{2}\left(\frac{y''}{y'}\right)^2 = 0.
$$ 

Given such a 3-parameter family of curves $\F$, one defines null cones whose rulings consist of the curves that are tangent to a fixed line at a fixed point. A smooth plane curve defines an associated curve in $\F$, and it is still true that this associated curve  is null. However, for a  general family $\F$, these cones may fail to be quadratic.

The families $\F$ for which the nullity condition is quadratic, and hence defines a conformal Lorentzian metric on $\F$,   are characterized  by a complicated non-linear PDE on the function $F(x,y,y',y'')$ that defines this family. This condition was
studied by K. W\"unschmann \cite{W}  and later  by \'E. Cartan \cite{Ca} and S.S. Chern \cite{Ch1, Ch2},  using the   method of equivalence. For  recent  presentations  of this deep result  see \cite{FKN,N}. 

\bigskip\centerline {*\qquad*\qquad*}
\medskip 

%\note{I added this paragraph (Connor's contribution). Please check!}
%
We end  with the  observation that  the  Tait-Knesser and  4 vertex theorems for Kepler conics 
(Theorems \ref{thm:TKKe} and \ref{thm:4v}) can be derived  from their Euclidean analogues via projective duality. Here is a sketch (more details  will appear in  \cite{BJ}). The equation $ax+by=1$ associates to  a point  in the $ab$ plane a line in the $xy$ plane and vice versa. To a   curve $C$ in one plane corresponds   its dual curve $C^*$ in the other plane, whose points parametrize the  lines tangent to $C$. One then shows that the dual of a Kepler conic is a circle and that duality preserves nesting and order of contacts of curves. It follows that the dual of the   osculating Kepler conic to $C$ is the  osculating circle to $C^*$ and the same holds for hyperosculating conics. Thus duality interchanges Euclidean and Keplerian vertices, reducing  the  Kepler version of the Tait-Knesser and  4 vertex theorems to their Euclidean analogues. 

\bigskip\noindent
{\bf Acknowledgements}.   After the first version of this article was posted, R. Pacheco and M. Salvai brought their relevant papers \cite{NP,Sa}
to our attention; we are grateful to them for it. We thank A. Izosimov for interesting discussions and the referee for  useful suggestions. 
GB was supported by CONACYT Grant A1-S-4588. ST was supported by NSF grant DMS-2005444.


\begin{thebibliography}{99}

\bibitem{Ar} V. I. Arnold. {\em Huygens and Barrow, Newton and Hooke: Pioneers in mathematical analysis and catastrophe theory from evolvents to quasicrystals.} Birkh\"auser Verlag, Basel, 1990. 

\bibitem{Be} J. Bertrand. {\it Th\'eor\`eme relatif au mouvement d'un point attir\'e vers un centre fixe}. C. R. Acad. Sci. {\bf 77} (1873), 849--853.

\bibitem{Bo} M. K. Bohlin. {\it Note sur le probl\`eme des deux corps et sur une int\'egration nouvelle dans le
probl\`eme des trois corps}. Bull. Astron. {\bf 28} (1911), 113--119.

\bibitem{BJ} G. Bor, C. Jackman. {\em Geometry and symmetries of Kepler orbits}. In preparation. 

\bibitem{Ca} \'E. Cartan. {\em La geometr\'ia de las ecuaciones diferenciales de tercer orden}. Rev. Mat. Hispano-Amer. {\bf 4} (1941), 1--31.

\bibitem{Ch1} S.-S. Chern. {\em Sur la g\'eom\'etrie d'une \'equation diff\'erentielle du troisi\`eme ordre}. C. R. Acad. Sci., Paris {\bf 204} (1937), 1227--1229.

\bibitem{Ch2} S.-S. Chern.{\em The geometry of the differential equations} $y''' = F(x,y,y',y'')$. Sci. Rep. Nat. Tsing Hua Univ.
{\bf 4} (1940), 97--111.

\bibitem{Giv} A. Givental. {\em  Kepler's Laws and Conic Sections}. Arnold Math J. {\bf 2} (2016), 139--148. 


\bibitem{GTT} E. Ghys, S. Tabachnikov, V. Timorin. {\em Osculating Curves: Around the Tait-Kneser Theorem}. Math. Intelligencer {\bf 35} (2013), no. 1, 61--66. 

\bibitem{Ka1} E. Kasner. {\it The trajectories of dynamics.} Trans. Amer. Math. Soc. {\bf 7} (1906), 401--424. 

\bibitem{Ka2} E. Kasner. {\it Differential-geometric aspects of dynamics}. The Princeton Colloquium, v. 3, AMS, Providence, 1913.

\bibitem{FKN}S. Frittelli, C. Kozameh, E.T. Newman. {\em Differential geometry from differential equations}, Comm. in Math. Phys. {\bf 223.2} (2001),  383--408.

\bibitem{NP} B. Nolasco, R. Pacheco. {\it Evolutes of plane curves and null curves in Minkowski 3-space.}
J. Geom. {\bf 108} (2017), 195--214. 

\bibitem{N}P. Nurowski. {\em Differential equations and conformal structures},  J. Geom. Phys. {\bf  55.1} (2005), 19--49.

\bibitem{OT} V. Ovsienko, S. Tabachnikov. {\it Projective differential geometry old and new. 
From the Schwarzian derivative to the cohomology of diffeomorphism groups. }
Cambridge Univ. Press, Cambridge, 2005.

\bibitem{Sa} M. Salvai. {\it Centro-affine invariants and the canonical Lorentz metric on the space of centered ellipses.} Kodai Math. J. {\bf 40} (2017), 21--30. 

\bibitem{T}P. G. Tait. {\em Note on the circles of curvature of a plane curve}. Proc. Edinburgh Math. Soc. {\bf 14} (1896), 403.

\bibitem{W}K. W\"unschmann. {\em \"Uber Ber\"uhrungsbedingungen bei Integralkurven von Differentialgleichungen}. Inauguraldissertation, Leipzig, Teubner, 1905, 6--13. 

\end{thebibliography}
\end{document}